\documentclass[10pt]{article}
\usepackage{amssymb}
\usepackage{amsmath}
\newtheorem{theorem}{Theorem}[section]
\newtheorem{proposition}[theorem]{Proposition}
\newtheorem{lemma}[theorem]{Lemma}
\newtheorem{corollary}[theorem]{Corollary}

\newtheorem{definition}[theorem]{Definition}
\newtheorem{remark}[theorem]{Remark}
\newtheorem{example}[theorem]{Example}

\newcommand{\bthm}{\begin{theorem}}
\newcommand{\ethm}{\end{theorem}}
\newcommand{\bpr}{\begin{proposition}}
\newcommand{\epr}{\end{proposition}}
\newcommand{\blem}{\begin{lemma}}
\newcommand{\elem}{\end{lemma}}
\newcommand{\bco}{\begin{corollary}}
\newcommand{\eco}{\end{corollary}}
\newcommand{\bde}{\begin{definition}\rm}
\newcommand{\ede}{\end{definition}}
\newcommand{\bre}{\begin{remark}\rm}
\newcommand{\ere}{\end{remark}}
\newcommand{\bex}{\begin{example}\rm}
\newcommand{\eex}{\end{example}}
\newcommand{\bprf}{\noindent{\it Proof.\ }}
\newcommand{\eprf}{\hspace*{\fill} \rule{1.6mm}{3.2mm} \vspace{1.6mm}}
\newcommand{\benu}{\begin{enumerate}\renewcommand{\labelenumi}{{\rm (\roman{enumi})}}\renewcommand{\itemsep}{0pt}}
\newcommand{\eenu}{\end{enumerate}}

\newcommand{\N}{\mathbb{N}}

\newcommand{\R}{\mathbb{R}}
\newcommand{\C}{\mathbb{C}}
\newcommand{\T}{\mathbb{T}}
\newcommand{\e}{\varepsilon}
\newcommand{\K}{\mathbb{K}}
\newcommand{\M}{\mathbb{M}}
\newcommand{\I}{\mathbb{I}}
\newcommand{\J}{\mathbb{J}}

\DeclareMathOperator{\spa}{span}
\DeclareMathOperator{\cspa}{\overline{span}}
\newcommand{\On}{{\cal O}_n}
\newcommand{\Oi}{{\cal O}_\infty}
\newcommand{\cp}{{\cal O}_n{\rtimes_{\alpha^\omega}}G}
\newcommand{\cpi}{{\cal O}_\infty{\rtimes_{\alpha^\omega}}G}
\newcommand{\cpr}{{\cal O}_n{\rtimes_{\alpha^\omega}}\mathbb{R}}

\newcommand{\ip}[2]{\langle\,{#1}\,|\,{#2}\,\rangle}
\newcommand{\W}{{\cal W}}

\begin{document}
\title{AF-embeddability of crossed products of Cuntz algebras}
\author{Takeshi KATSURA \\
Department of Mathematical Sciences\\
University of Tokyo, Komaba, Tokyo, 153-8914, JAPAN\\
e-mail: {\tt katsu@ms.u-tokyo.ac.jp}}
\date{}

\maketitle

\begin{abstract}
{\footnotesize 
We investigate crossed products of Cuntz algebras by quasi-free actions 
of abelian groups. 
We prove that our algebras are AF-embeddable 
when actions satisfy a certain condition.
We also give a necessary and sufficient condition that our algebras become 
simple and purely infinite, 
and consequently our algebras are either purely infinite or AF-embeddable 
when they are simple.}
\end{abstract}

\section{Introduction}\label{INTRO}

There had been no examples of simple $C^*$-algebras which have both 
a finite projection and an infinite one 
until M. R\o rdam found such a $C^*$-algebra recently \cite{Ro}.
However, we have found no examples of such simple $C^*$-algebras 
among nuclear ones, so far.
Moreover we have not known examples of simple nuclear $C^*$-algebras
which are not stably finite nor purely infinite.
The property `stable finiteness' has recently attracted much attention
in connection with quasidiagonality and AF-embeddability.
It is easy to see that AF-embeddability implies quasidiagonality 
and that quasidiagonality implies stable finiteness.
It is still open whether or not stable finiteness implies AF-embeddability 
for nuclear $C^*$-algebras.
On this topic, there is a nice survey \cite{Br3} written by N. P. Brown.
Since M. Pimsner and D. Voiculescu showed 
that the irrational rotation algebras are AF-embeddable \cite{PV}, 
several authors have studied AF-embeddability 
of some classes of $C^*$-algebras. 
In particular, we can find many papers dealing with AF-embeddability of 
crossed products of finite $C^*$-algebras, for example, 
\cite{Pu}, \cite{Pi1}, \cite{Pi2} for those of commutative $C^*$-algebras, 
and \cite{V}, \cite{Br1}, \cite{Br2} for those of AF-algebras.
On the other hand, the author has been unable to find any article 
related to AF-embeddability of crossed products of infinite $C^*$-algebras.
We remark that it seems more difficult to show AF-embeddability 
of crossed products of infinite $C^*$-algebras by continuous groups 
than those of finite $C^*$-algebras.
For crossed products of finite $C^*$-algebras, 
there is a method to derive AF-embeddability of crossed products 
by continuous groups from the discrete group case 
by using Green's imprimitivity theorem (\cite{Gr}, see also \cite{Br2}).
However, for infinite $C^*$-algebras, we cannot use this method 
because their crossed products by discrete groups are never embedded 
into AF-algebras.

In this paper, we will deal with crossed products of Cuntz algebras 
by quasi-free actions of abelian groups, 
whose ideal structures were examined in our previous paper \cite{Ka}.
We will prove the AF-embeddability of our algebras under a certain
condition for actions.
To the author's knowledge, this is the first case to have succeeded in 
embedding crossed products of purely infinite $C^*$-algebras 
into AF-algebras except trivial cases.
We will also show that our algebras are either purely infinite 
or AF-embeddable when they are simple.

This paper is organized as follows.
After some preliminaries, 
we will show that the crossed products are AF-embeddable 
when actions satisfy a certain condition (Theorem \ref{embed}).
They were known to be stably finite in the case that the group is the 
real number group $\R$ \cite{KK1}.
In the case that the group is compact, this condition is also sufficient
for the crossed products to be AF-embeddable, 
and moreover the crossed products become AF-algebras 
under this condition. 
For the general setting, we do not know whether our algebra is AF-embeddable 
or not when the action does not satisfy the condition 
(see Remark \ref{remark}).
In section \ref{secpi}, 
we will give a necessary and sufficient condition that our algebras become 
simple and purely infinite. 
Combining this characterization with our result on AF-embeddability, 
we can easily get the dichotomy which says that our algebras are 
either purely infinite or AF-embeddable when they are simple.
In the last section, we will deal with crossed products of the Cuntz algebra
$\Oi$, which is generated by infinitely many isometries, by the same type of 
actions of abelian groups.
We will prove AF-embeddability of such algebras under a certain condition 
for actions, 
and give a necessary and sufficient condition for such algebras to be 
simple and purely infinite which will be shown to be equivalent to the
property that they are simple.

\medskip

{\bf Acknowledgment.} This work was done while the author stayed in the Mathematical Sciences Research Institute, and he would like to express his gratitude to MSRI for their hospitality.
The author is grateful to Masaki Izumi for various comments and many important suggestions, and to his advisor Yasuyuki Kawahigashi for his support and encouragement.
This work was partially supported by Research Fellowship for Young Scientists of the Japan Society for the Promotion of Science (JSPS) and Honda Heizaemon memorial fellowship of the Japan Association for Mathematical Sciences.

\section{Preliminaries}\label{PRE}

In this section, we review some results and fix the notation.
For $n=2,3,\ldots$, the Cuntz algebra $\On$ is the universal $C^*$-algebra generated by $n$ isometries $S_1,S_2,\ldots,S_n$, satisfying $\sum_{i=1}^nS_iS_i^*=1$.
For $k\in\N=\{0,1,\ldots\}$, we define the set $\W_n^{(k)}$ of $k$-tuples 
by $\W_n^{(0)}=\{\emptyset\}$ and
$$\W_n^{(k)}=\big\{ (i_1,i_2,\ldots,i_k)\ \big|\ i_j\in\{1,2,\ldots,n\}\big\}.$$
We set $\W_n=\bigcup_{k=0}^\infty\W_n^{(k)}$. 
For $\mu=(i_1,i_2,\ldots,i_k)\in\W_n$, we denote its length $k$ by $|\mu|$, and set $S_\mu=S_{i_1}S_{i_2}\cdots S_{i_k}\in\On$. 
Note that $|\emptyset|=0,\ S_\emptyset=1$. 
For $\mu=(i_1,i_2,\ldots,i_k),\nu=(j_1,j_2,\ldots,j_l)\in\W_n$, we define their product $\mu\nu\in\W_n$ by $\mu\nu=(i_1,i_2,\ldots,i_k,j_1,j_2,\ldots,j_l)$.

We fix a locally compact abelian group $G$ whose dual group is denoted by 
$\Gamma$ which is also a locally compact abelian group.
We always use $+$ for multiplicative operations of abelian groups except for 
$\T$, which is the group of the unit circle in the complex plane $\C$. 
The pairing of $t\in G$ and $\gamma\in\Gamma$ is denoted by 
$\ip{t}{\gamma}\in\T$. 

\bde\label{wga}
Let $\omega=(\omega_1,\omega_2,\ldots,\omega_n)\in\Gamma^n$ be given. 
We define the action $\alpha^\omega:G\curvearrowright\On$ by
$$\alpha^\omega_t(S_i)=\ip{t}{\omega_i}S_i\quad (i=1,2,\ldots,n,\ t\in G).$$
\ede

This type of action is called quasi-free 
(see \cite{E} for quasi-free actions on the Cuntz algebras).
Since the abelian group $G$ is amenable, the reduced crossed product of 
the action $\alpha^\omega:G\curvearrowright\On$ coincides with 
the full crossed product of it.
We will denote it by $\cp$ and call it the crossed product.
The crossed product $\cp$ has a $C^*$-subalgebra 
$\C 1{\rtimes_{\alpha^\omega}}G$, which is isomorphic to $C_0(\Gamma)$.
Throughout this paper, 
we always consider $C_0(\Gamma)$ as a $C^*$-subalgebra of $\cp$, 
and use $f,g,\ldots$ for denoting elements of $C_0(\Gamma)\subset \cp$.
The Cuntz algebra $\On$ is naturally embedded 
into the multiplier algebra $M(\cp)$ of $\cp$. 
For each $\mu=(i_1,i_2,\ldots,i_k)$ in $\W_n$, 
we define an element $\omega_\mu$ of $\Gamma$ by 
$\omega_\mu=\sum_{j=1}^{k}\omega_{i_j}$.
For $\gamma_0\in\Gamma$, we define a (reverse) shift automorphism 
$\sigma_{\gamma_0}:C_0(\Gamma)\to C_0(\Gamma)$ by 
$(\sigma_{\gamma_0} f)(\gamma)=f(\gamma+\gamma_0)$ for $f\in C_0(\Gamma)$. 
Once noting that $\alpha^\omega_t(S_\mu)=\ip{t}{\omega_\mu}S_\mu$ for $\mu\in\W_n$, one can easily verify that $fS_\mu =S_\mu\sigma_{\omega_\mu}f$ for any $f\in C_0(\Gamma)\subset \cp$ and any $\mu\in\W_n$.
The linear span of $\{ S_\mu fS_\nu^*\mid \mu,\nu\in\W_n,\ f\in C_0(\Gamma)\}$ is dense in $\cp$ (see \cite{Ka}).
We denote by $\M_k$ the $C^*$-algebra of $k \times k$ matrices for $k=1,2,\ldots$, and by $\K$ the $C^*$-algebra of compact operators of the infinite dimensional separable Hilbert space.

\section{AF-embeddability of $\cp$}\label{AFE}

A. Kishimoto and A. Kumjian showed that $\cpr$ is stably projectionless if all the $\omega_i$'s have the same sign by using the KMS-state \cite[Theorem 4.1]{KK1}.
Thus $\cpr$ is stably finite in this case.
In this section, we will show that $\cp$ becomes AF-embeddable if $\omega$ satisfies a certain condition.
This gives another proof of the stable finiteness of $\cpr$ when all the $\omega_i$'s have the same sign.
More precisely, we will prove that if $-\omega_i\notin\overline{\{\omega_{\mu}\mid\mu\in\W_n\}}$ for any $i\in\{1,2,\ldots,n\}$, then $\cp$ is AF-embeddable (Theorem \ref{embed}).
Here we note that $\overline{\{\omega_{\mu}\mid\mu\in\W_n\}}$ is 
the closed semigroup generated by $\omega_1,\omega_2,\ldots,\omega_n$.

Let us take a faithful representation $\On\hookrightarrow B(H)$ 
for some Hilbert space $H$.
There exists a canonical embedding $\cp\hookrightarrow B(H\otimes L^2(G))$.
Since $L^2(G)$ is isomorphic to $L^2(\Gamma)$ via the Fourier transform, 
we can consider $\cp$ as a subalgebra of $B(H\otimes L^2(\Gamma))$.
In this setting, an element of $C_0(\Gamma)\subset\cp$ acts by multiplication 
on $L^2(\Gamma)$ and as identity on $H$.
Note that the weak closure of $C_0(\Gamma)$ in $B(H\otimes L^2(\Gamma))$ 
is $L^\infty(\Gamma)$.

Throughout this section, we fix $\omega\in\Gamma^n$ satisfying 
$-\omega_i\notin\overline{\{\omega_{\mu}\mid\mu\in\W_n\}}$ for any $i$.
We also fix an open base $\{U_i\}_{i\in\I}$ of $\Gamma$ 
such that for any $i\in\I$, $\overline{U_i}$ is compact and for any 
$i\in\I$ and $\mu\in\W_n$, there exists $j\in\I$ with 
$U_j=U_i-\omega_\mu$.
Obviously such an open base exists, and we can take countable one 
when $\Gamma$ satisfies the second countability axiom.
For each $i\in\I$, 
let us consider the characteristic function $\chi_{U_i}$ of $U_i$ 
which is an element of $L^{\infty}(\Gamma)\subset B(H\otimes L^2(\Gamma))$.
Let $D_0(\Gamma)$ be the $C^*$-algebra generated by $\{\chi_{U_i}\}_{i\in\I}$. 
Let us denote by $\Lambda$ the directed set of all finite subsets of $\I$ 
whose order is defined by the inclusion. 
For $\lambda=\{i_1,i_2,\ldots,i_k\}\in\Lambda$, the $C^*$-subalgebra 
$D_{\lambda}$ of $D_0(\Gamma)$ is defined by the $C^*$-algebra generated by 
$\chi_{U_{i_1}},\chi_{U_{i_2}},\ldots,\chi_{U_{i_k}}$.
One can easily verify the following.

\blem\label{D_0}
\benu
\item $C_0(\Gamma)\subset D_0(\Gamma)$.
\item We can define the shift $*$-homomorphism 
$\sigma_{\omega_{\mu}}:D_0(\Gamma)\to D_0(\Gamma)$ for any $\mu\in\W_n$.
\item $\varinjlim D_{\lambda}=D_0(\Gamma)$.
\item $D_0(\Gamma)$ is an AF-algebra.
\eenu
\elem

Define a subspace $A$ of $B(H\otimes L^2(\Gamma))$ by
$$A=\cspa\{S_\mu fS_\nu^*\mid \mu,\nu\in\W_n,\ f\in D_0(\Gamma)\}.$$
By Lemma \ref{D_0} (ii), $A$ is a $C^*$-algebra and 
by Lemma \ref{D_0} (i), $A$ contains $\cp$.
We will show that $A$ is an AF-algebra when 
$-\omega_i\notin\overline{\{\omega_{\mu}\mid\mu\in\W_n\}}$ 
for any $i\in\{1,2,\ldots,n\}$, 
which implies that $\cp$ is AF-embeddable.
We denote by $A_\lambda$ the $C^*$-algebra generated by 
$\{S_\mu \chi_{U_i}S_\nu^*\mid \mu,\nu\in\W_n,\ i\in\lambda\}$.
It is easy to see the following.

\blem\label{indlim}
With the above notation, we have $A=\varinjlim A_\lambda$.
\elem

By Lemma \ref{indlim}, to prove that $A$ is an AF-algebra, 
it suffices to show that $A_\lambda$ is an AF-algebra 
for any $\lambda\in\Lambda$.
Let us take $\lambda\in\Lambda$ arbitrarily, and fix it.
Let $p_1,p_2,\ldots,p_L$ be minimal projections of $D_\lambda$
and $p=\sum_{l=1}^L p_l$ be its unit.
Note that $A_\lambda$ is generated by $\{S_\mu p_l S_\nu^*\mid \mu,\nu\in\W_n,\ l=1,2,\ldots,L\}$.
Only in the next lemma, we use directly the assumption that $\omega$ satisfies 
$-\omega_i\notin\overline{\{\omega_{\mu}\mid\mu\in\W_n\}}$ for any 
$i\in\{1,2,\ldots,n\}$, and this lemma implies all the following lemmas and 
the fact that $A_\lambda$ is an AF-algebra.

\blem\label{shift}
There exists $K\in\N$ such that $pS_\mu p=0$ 
for any $\mu\in\W_n$ with $|\mu|>K$.
\elem

\bprf
If we define a subset $U=\bigcup_{i\in\lambda}U_i$ of $\Gamma$,
then $p$ is the characteristic function of $U$.
The closure of $U$ is compact since $\overline{U_i}$ is compact 
for any $i\in\lambda$.
To derive a contradiction, assume that for any $k\in\N$, 
there exists $\mu_k\in\W_n$ such that $|\mu_k|>k$ and $pS_{\mu_k} p\neq 0$.
Then we have $S_{\mu_k}^*pS_{\mu_k} p\neq 0$.
Since $S_{\mu_k}^*pS_{\mu_k}$ is the characteristic function of 
$U-\omega_{\mu_k}$, 
there exists $\gamma_k\in (U-\omega_{\mu_k})\cap U$.
We have $\omega_{\mu_k}=(\gamma_k+\omega_{\mu_k})-\gamma_k\in U-U$ 
for any $k\in\N$.
Since $\overline{U-U}$ is compact, there exists an increasing subsequence 
$k_1,k_2,\ldots,k_m,\ldots$ of $\N$ such that $\omega_{\mu_{k_m}}$ converges 
to some element $\gamma_0\in\overline{U-U}$ when $m$ goes to infinity.
By replacing it by a subsequence of $\{k_m\}$ if necessary,
we may assume that the number of $i$ appearing in $\mu_{k_m}\in\W_n$ 
does not decrease for $i=1,2,\ldots,n$.
Since $|\mu_{k_m}|\to\infty$ when $m\to\infty$, there exists 
$i_0\in\{1,2,\ldots,n\}$ such that the number of $i_0$ appearing 
in $\mu_{k_m}$ diverges to infinity when $m\to\infty$.
By replacing it by a subsequence of $\{k_m\}$ if necessary, 
we may assume that the number of $i_0$ appearing in $\mu_{k_m}$ 
increases strictly.
Thus, we have $\omega_{\mu_{k_m}}-\omega_{\mu_{k_{m-1}}}-\omega_{i_0}
\in\{\omega_{\mu}\mid\mu\in\W_n\}$ for any $m\in\N$.
By
$$\lim_{m\to\infty}(\omega_{\mu_{k_m}}-\omega_{\mu_{k_{m-1}}}-\omega_{i_0})=\gamma_0-\gamma_0-\omega_{i_0}=-\omega_{i_0},$$
we have $-\omega_{i_0}\in\overline{\{\omega_{\mu}\mid\mu\in\W_n\}}$.
This is a contradiction.
\eprf

We fix a positive integer $K$ satisfying the condition in Lemma \ref{shift}.
Before going further, 
we remark that $S_\mu fS_\mu^*$ and $S_\nu gS_\nu^*$ commute 
for any $\mu,\nu\in\W_n$ and any $f,g\in D_0(\Gamma)$.
This fact will be used without further notice.
For $k\in\N$, define a $*$-endomorphism $\rho_k$ of $A_\lambda$ by 
$\rho_k(x)=\sum_{\mu\in\W_n^{(k)}}S_\mu xS_\mu^*$.
Set a projection $q$ in $A_\lambda$ by
$$q=\left(\prod_{k=1}^K\big(1-\rho_k(p)\big)\right)p.$$
Note that $q\leq p$ and $q\leq 1-\rho_k(p)$ for $k=1,2,\ldots,K$.

\blem\label{divide1}
For any $\mu,\nu\in\W_n$ with $\mu\neq\nu$, two projections $S_\mu qS_\mu^*$ and $S_\nu qS_\nu^*$ are orthogonal to each other.
\elem

\bprf
It suffices to show that $qS_\mu q=0$ for any $\mu\in\W_n$ with $\mu\neq\emptyset$.
When $1\leq |\mu|\leq K$, since
$$\left(1-\rho_{|\mu|}(p)\right)S_\mu p=(S_\mu-S_\mu p)p=0,$$
we have $qS_\mu q=0$.
When $|\mu|>K$, $pS_\mu p=0$ by Lemma \ref{shift}, so $q S_\mu q=0$. 
\eprf

Denote a set $\{\mu\in\W_n\mid |\mu|\leq K\}$ by $\W$.

\blem\label{divide2}
We have $\sum_{\mu\in\W}S_\mu qS_\mu^*p=p$.
\elem

\bprf
For $l=1,2,\ldots,K$, we have
$$\rho_l(q)=\left(\prod_{k=1}^K\big(1-\rho_{l+k}(p)\big)\right)\rho_l(p).$$
Since $\left(1-\rho_k(p)\right)p=p$ for $k>K$ by Lemma \ref{shift},
we have 
$$\rho_l(q)p=\left(\prod_{k=l+1}^K\big(1-\rho_{k}(p)\big)\right)\rho_l(p)p.$$
Hence 
$$\sum_{\mu\in\W}S_\mu qS_\mu^*p=\sum_{l=0}^K\sum_{\mu\in\W_n^{(l)}}\rho_l(q)p
=\sum_{l=0}^K\left(\prod_{k=l+1}^K\big(1-\rho_{k}(p)\big)\right)\rho_l(p)p=p.$$
\eprf

Let us define a projection $p_0$ by $p_0=1-p=1-\sum_{l=1}^L p_l$ 
where $p_1,p_2,\ldots,p_L$ are the minimal projections of $D_\lambda$.
Note that $p_0,p_1,\ldots,p_L$ is a set of mutually orthogonal projections 
whose sum is $1$.
Let $\J'$ be a set of all maps from the set 
$\W=\{\mu\in\W_n\mid |\mu|\leq K\}$ to the set $\{0,1,2,\ldots,L\}$.
For $\tau\in\J'$, we define a projection $q_{\tau}\in A_\lambda$ by
$$q_{\tau}=q\prod_{\mu\in\W}S_\mu^*p_{\tau(\mu)}S_\mu.$$
Set $\J=\{\tau\in\J'\mid q_{\tau}\neq 0\}$.

\blem\label{divide3}
\benu
\item $\{q_{\tau}\}_{\tau\in\J}$ is a set of mutually orthogonal non-zero projections.
\item $\sum_{\tau\in\J}q_{\tau}=q$.
\item For $\mu\in\W$, $\tau\in\J$ and $l\in\{1,2,\ldots,L\}$, we have $S_\mu q_{\tau} S_\mu^*p_l=\delta_{\tau(\mu),l}S_\mu q_{\tau} S_\mu^*$.
\eenu
\elem

\bprf
\benu
\item If $\tau_1\neq\tau_2$, then $\tau_1(\mu)\neq\tau_2(\mu)$ for some 
$\mu\in\W$. 
Now $q_{\tau_1} q_{\tau_2}=0$ follows from
$$(S_\mu^* p_{\tau_1(\mu)}S_\mu)(S_\mu^* p_{\tau_2(\mu)}S_\mu)=S_\mu^*S_\mu S_\mu^* p_{\tau_1(\mu)}p_{\tau_2(\mu)}S_\mu=0,$$
since $q_{\tau_1}\leq S_\mu^* p_{\tau_1(\mu)}S_\mu$ and $q_{\tau_2}\leq S_\mu^* p_{\tau_2(\mu)}S_\mu$.
\item $\sum_{\tau\in\J}q_{\tau}=\sum_{\tau\in\J'}q_{\tau}=q\prod_{\mu\in\W}S_\mu^*(p_0+p_1+\cdots+p_L)S_\mu=q$.
\item It follows from the fact that 
$$S_\mu (S_\mu^* p_{\tau(\mu)} S_\mu) S_\mu^*p_l=S_\mu S_\mu^* p_{\tau(\mu)}p_l S_\mu S_\mu^*=\delta_{\tau(\mu),l}S_\mu (S_\mu^* p_{\tau(\mu)} S_\mu) S_\mu^*.$$
\eprf
\eenu

\bpr\label{K}
We have $A_\lambda\cong\bigoplus_{\tau\in\J}\K$. Hence, $A_\lambda$ is an AF-algebra.
\epr

\bprf
For any $\tau_1,\tau_2\in\J$ and $\mu,\nu\in\W_n$ with $\mu\neq\nu$, 
we have $(S_\mu q_{\tau_1} S_\mu^*)(S_\nu q_{\tau_2} S_\nu^*)=0$ 
by Lemma \ref{divide1}.
Thus, for $\tau_1,\tau_2\in\J$ and $\mu_1,\mu_2,\nu_1,\nu_2\in\W_n$, we get
\begin{align*}
(S_{\mu_1}q_{\tau_1} S_{\nu_1}^*)(S_{\mu_2}q_{\tau_2} S_{\nu_2}^*)&=\delta_{\nu_1,\mu_2}S_{\mu_1}q_{\tau_1} q_{\tau_2} S_{\nu_2}^*\\
&=\delta_{\nu_1,\mu_2}\delta_{\tau_1,\tau_2}S_{\mu_1}q_{\tau_1} S_{\nu_2}^*.
\end{align*}
For any $\tau\in\J$, the set $\{S_{\mu}q_{\tau} S_{\nu}^*\}_{\mu,\nu\in\W_n}$ 
satisfies the relation of matrix units, so the $C^*$-algebra generated by 
$\{S_{\mu}q_{\tau} S_{\nu}^*\}_{\mu,\nu\in\W_n}$ is isomorphic to $\K$.
For any two elements $\tau_1,\tau_2$ in $\J$, the $C^*$-algebra generated by 
$\{S_{\mu}q_{\tau_1} S_{\nu}^*\}_{\mu,\nu\in\W_n}$ is orthogonal to the 
$C^*$-algebra generated by $\{S_{\mu}q_{\tau_2} S_{\nu}^*\}_{\mu,\nu\in\W_n}$.
Therefore, the $C^*$-algebra generated by 
$\{S_{\mu}q_{\tau} S_{\nu}^*\mid \mu,\nu\in\W_n,\ \tau\in\J\}$ 
is isomorphic to $\bigoplus_{\tau\in\J}\K$.

Since $q_{\tau}\in A_\lambda$ for any $\tau\in\J$, 
the $C^*$-algebra generated by 
$\{S_{\mu}q_{\tau} S_{\nu}^*\mid \mu,\nu\in\W_n,\ \tau\in\J\}$ 
is contained in $A_\lambda$.
Conversely, for $l=1,2,\ldots,L$, 
\begin{align*}
p_l&=pp_l\\
&=\sum_{\mu\in\W}S_\mu qS_\mu^*pp_l&&(\mbox{by Lemma \ref{divide2}})\\
&=\sum_{\mu\in\W}S_\mu \left(\sum_{\tau\in\J}q_{\tau}\right)S_\mu^* p_l&&(\mbox{by Lemma \ref{divide3} (ii)})\\
&=\sum_{\mu\in\W,\tau\in\J}S_\mu q_{\tau}S_\mu^*p_l&\\
&=\sum_{\stackrel{\mbox{$\scriptstyle \mu\in\W,\tau\in\J$,}}{\mbox{$\scriptstyle \mbox{s.t. }\tau(\mu)=l$}}}S_\mu q_{\tau}S_\mu^*&&(\mbox{by Lemma \ref{divide3} (iii)}).
\end{align*}
Thus, for any $\mu,\nu\in\W_n$ and $l=1,2,\ldots,L$, 
the element $S_\mu p_lS_\nu^*$ is contained in the $C^*$-algebra generated by 
$\{S_{\mu}q_{\tau} S_{\nu}^*\mid \mu,\nu\in\W_n,\ \tau\in\J\}$.
Therefore $A_\gamma$ coincides with the $C^*$-algebra generated by 
$\{S_{\mu}q_\tau S_{\nu}^*\mid \mu,\nu\in\W_n,\ \tau\in\J\}$ 
which was proved to be isomorphic to $\bigoplus_{\tau\in\J}\K$.
\eprf

Now we can prove the main theorem.

\bthm\label{embed}
If $\omega$ satisfies 
$-\omega_i\notin\overline{\{\omega_{\mu}\mid\mu\in\W_n\}}$ 
for any $i\in\{1,2,\ldots,n\}$, 
then the crossed product $\cp$ is AF-embeddable.
\ethm

\bprf
The $C^*$-algebra $A$ is an AF-algebra 
because it is an inductive limit of AF-algebras.
Since the crossed product $\cp$ is naturally embedded into $A$, 
it is AF-embeddable.
\eprf

\bpr\label{AF!}
When $G$ is compact, the following are equivalent:
\benu
\item $-\omega_i\notin\overline{\{\omega_{\mu}\mid\mu\in\W_n\}}$ for any $i\in\{1,2,\ldots,n\}$.
\item The crossed product $\cp$ is stably finite.
\item The crossed product $\cp$ is AF-embeddable.
\item The crossed product $\cp$ itself is an AF-algebra.
\eenu
\epr

\bprf
(i)$\Rightarrow$(iv): 
Note that $\Gamma$ is discrete when $G$ is compact.
We can take $\big\{\{\gamma\}\big\}_{\gamma\in\Gamma}$ 
for an open base $\{U_i\}_{i\in\I}$.
Then the $C^*$-algebra $A$ which was proved to be an AF-algebra 
is $\cp$ itself.
Thus, $\cp$ is an AF-algebra.

(iv)$\Rightarrow$(iii)$\Rightarrow$(ii): Obvious.

(ii)$\Rightarrow$(i): 
If there exists $i\in\{1,2,\ldots,n\}$ such that 
$-\omega_i\in\overline{\{\omega_{\mu}\mid\mu\in\W_n\}}$, 
then there exists $\mu'\in\W_n$ with $-\omega_i=\omega_{\mu'}$.
Hence $\mu=i\mu'\in\W_n$ satisfies $|\mu|\geq 1$ and $\omega_\mu=0$.
Set $u=S_\mu\chi\in\cp$ where $\chi\in C_0(\Gamma)$ is the characteristic 
function of $\{0\}$.
We have $u^*u=\chi$ and $uu^*=S_\mu\chi S_\mu^*$.
We get $\chi\neq S_\mu\chi S_\mu^*$ from $|\mu|\geq 1$, 
and $\chi (S_\mu\chi S_\mu^*)=S_\mu\chi S_\mu^*$ from $\omega_\mu=0$.
Therefore $\chi$ is an infinite projection.
Thus $\cp$ is not stably finite.
\eprf

\bre\label{remark}
When $G=\R$, Theorem \ref{embed} implies that $\cpr$ is AF-embeddable
if all the $\omega_i$'s have the same sign.
If there exist $i,j$ such that $\omega_i<0<\omega_j$,
then $\cpr$ has infinite projections hence it is not AF-embeddable.
We do not know whether $\cpr$ is AF-embeddable or not if there exists 
$i\in\{1,2,\ldots,n\}$ such that $\omega_i=0$ and all the other 
$\omega_i$'s have the same sign,
though it is not hard to see that it is stably finite.
\ere

\section{Pure infiniteness of $\cp$}\label{secpi}

In this section, we investigate for which $\omega\in\Gamma^n$ 
the crossed product $\cp$ becomes simple and purely infinite.
Recall that a simple $C^*$-algebra is called purely infinite 
if any non-zero hereditary subalgebra has an infinite projection.
An element $x$ of a $C^*$-algebra is called a scaling element if 
$(x^*x)(xx^*)=xx^*$ and $x^*x\neq xx^*$.
In \cite{BC}, B. E. Blackadar and J. Cuntz showed that 
if a simple stable $C^*$-algebra has a scaling element, 
then it has an infinite projection.
One can omit the assumption of stability (Proposition \ref{bc1}).
To do so, we need the following standard lemma.

\blem\label{d}
Let $A$ be a $C^*$-algebra, $p$ a projection of $A$, and 
$a$ a positive element of $A$.
If there exist $x_1,x_2\ldots,x_K$ and $y_1,y_2,\ldots,y_K$ in $A$ with
$$\left\|p-\sum_{k=1}^K x_k a y_k\right\|<\frac12,$$
then there exist $z_1,z_2\ldots,z_{2K}$ in $A$ such that
$$p=\sum_{k=1}^{2K}z_k^* a z_k.$$

In particular, if $A$ is simple $C^*$-algebra, $p$ is a projection of $A$, 
and $a$ is a non-zero positive element of $A$, 
then there exist $x_1,x_2\ldots,x_{K}$ in $A$ such that 
$p=\sum_{k=1}^{K}x_k^* a x_k$.
\elem

\bprf
See \cite[Lemma V.5.4]{D}, for example.
\eprf

\bpr\label{bc1}
If a $C^*$-algebra $A$ is simple and has a scaling element, 
then it has an infinite projection.
\epr

\bprf
If $A$ has a scaling element,
then $A$ has mutually orthogonal, mutually equivalent, non-zero projections 
$\{p_k\}_{k=1}^\infty$ and a positive element $a$ with $ap_k=p_k$ for any $k$ 
\cite[Theorem 3.1]{BC}.
Since $A$ is simple, there exist $x_1,x_2\ldots,x_K$ and $y_1,y_2,\ldots,y_K$ 
in $A$ with
$$\left\|a-\sum_{k=1}^K x_k p_1 y_k\right\|<\frac12.$$
Let us set $p=\sum_{k=1}^{2K+1}p_k$, which is a projection.
Then we have 
$$\left\|p-\sum_{k=1}^K x_k p_1 (y_kp)\right\|
=\left\|\big(a-\sum_{k=1}^K x_k p_1 y_k\big)p\right\|<\frac12,$$
since $ap=p$.
Hence there exist $z_1,z_2,\ldots, z_{2K}$ in $A$ such that 
$p=\sum_{k=1}^{2K} z_k^* p_1 z_k$ by Lemma \ref{d}.
For $k=1,2,\ldots,2K$, let $u_k$ be a partial isometry 
with $u_k^*u_k=p_1, u_ku_k^*=p_k$.
Set $z=\sum_{k=1}^{2K}u_kz_k$.
Then we have $z^*z=\sum_{k=1}^{2K}z_k^* p_1 z_k=p$.
Since $zz^*(\sum_{k=1}^{2K}p_k)=zz^*$, 
we have $zz^*\leq \sum_{k=1}^{2K}p_k<p$.
Therefore $p$ is an infinite projection.
\eprf

A. Kishimoto and A. Kumjian proved that $\cpr$ is simple and 
purely infinite if and only if the closed semigroup generated by 
$\omega_1,\omega_2,\ldots,\omega_n$ is $\R$ in \cite{KK2}.
We will generalize their result for our setting 
by using the same technique as in \cite{KK2}.
Namely, we will prove that $\cp$ is simple and purely infinite 
if and only if the closed semigroup generated by 
$\omega_1,\omega_2,\ldots,\omega_n$ is $\Gamma$.
When $\omega$ satisfies $\Gamma=\overline{\{\omega_{\mu}\mid\mu\in\W_n\}}$, 
the crossed product $\cp$ is simple by \cite[Theorem 4.4]{Ki} 
(see also \cite[Theorem 4.8]{Ka}).
First we will show that $\cp$ has a scaling element 
and hence an infinite projection.

\blem\label{type-}
Suppose that $\omega$ satisfies 
$\Gamma=\overline{\{\omega_{\mu}\mid\mu\in\W_n\}}$.
For any neighborhood $U$ of $0\in\Gamma$ and any positive integer $K$,
there exist $K$ elements $\mu_1,\mu_2,\cdots,\mu_K$ of $\W_n$ such that $\omega_{\mu_k}\in U$ for $k=1,2,\ldots,K$ and $S_{\mu_{k}}^*S_{\mu_{l}}=\delta_{k,l}$.
\elem

\bprf
We can find $K$ elements $\nu_1,\nu_2,\ldots,\nu_K$ of $\W_n$ such that $S_{\nu_k}^*S_{\nu_l}=\delta_{k,l}$.
For $k=1,2,\ldots,K$, there exists $\nu_k'\in\W_n$ with $\omega_{\nu_k'}\in U-\omega_{\nu_k}$ because $U-\omega_{\nu_k}$ is open and $\{\omega_{\mu}\mid\mu\in\W_n\}$ is dense in $\Gamma$.
Set $\mu_k=\nu_k\nu_k'$ for $k=1,2,\ldots,K$.
Then $S_{\mu_k}^*S_{\mu_l}=\delta_{k,l}$ and $\omega_{\mu_k}=\omega_{\nu_k}+\omega_{\nu_k'}\in U$ for $k=1,2,\ldots,K$.
\eprf

\blem\label{partition}
Suppose that $\omega$ satisfies 
$\Gamma=\overline{\{\omega_{\mu}\mid\mu\in\W_n\}}$.
Let $X$ be a compact neighborhood of $0\in\Gamma$ that differs from $\Gamma$. 
Then, there exist positive functions $f_1,f_2,\ldots,f_K\in C_0(\Gamma)$ and $\mu_1,\mu_2,\ldots,\mu_K\in\W_n$ satisfying the following conditions:
\benu
\item $S_{\mu_k}^*S_{\mu_l}=\delta_{k,l}$.
\item $\sum_{k=1}^Kf_k(\gamma)=1$ for any $\gamma\in X$.
\item $\sum_{k=1}^Kf_k(\gamma_0)\neq 0,1$ for some $\gamma_0\in\Gamma$.
\item The support of $\sigma_{-\omega_{\mu_k}}f_k$ is contained in $X$ for $k=1,2,\ldots,K$.
\eenu
\elem

\bprf
Let us choose an open neighborhood $U_1$ of $0$ such that the open neighborhood $U=U_1+U_1$ of 0 is contained in $X$,
and then choose an open neighborhood $U_2$ of $0$ such that $\overline{U_2}\subset U_1$.
For any $\gamma\in\Gamma$, there exists $\mu\in\W_n$ with $\omega_\mu\in U_2+\gamma$ because $\{\omega_{\mu}\mid\mu\in\W_n\}$ is dense in $\Gamma$.
Therefore $\bigcup_{\mu\in\W_n}(U_2-\omega_\mu)=\Gamma$.
Since $X$ is compact, there exist finite elements $\nu_1,\nu_2,\ldots,\nu_K$ of $\W_n$ such that 
$$X\subsetneqq\bigcup_{k=1}^K (U_2-\omega_{\nu_k}).$$
By Lemma \ref{type-}, there exist $K$ elements $\nu_1',\nu_2',\ldots,\nu_K'\in \W_n$ such that $S_{\nu_k'}^*S_{\nu_l'}=\delta_{k,l}$ and $\omega_{\nu_k'}\in U_1$ for $k=1,2,\ldots,K$.
Set $\mu_k=\nu_k'\nu_k$ for $k=1,2,\ldots,K$.
Then $S_{\mu_k}^*S_{\mu_l}=\delta_{k,l}$.
For $k=1,2,\ldots,K$, we get
\begin{align*}
\overline{U_2-\omega_{\nu_k}}&\subset U_1-\omega_{\nu_k}\\
&\subset U_1+U_1-\omega_{\nu_k}-\omega_{\nu_k'}\\
&=U-\omega_{\mu_k},
\end{align*}
since $\overline{U_2}\subset U_1$ and $\omega_{\nu_k'}\in U_1$.
For $k=1,2,\ldots,K$, let $g_k\in C_0(\Gamma)$ be a function with $0\leq g_k\leq 1$ such that $g_k(\gamma)=1$ for $\gamma\in \overline{U_2-\omega_{\nu_k}}$ and $g_k(\gamma)=0$ for $\gamma\notin U-\omega_{\mu_k}$.
Let us choose a continuous positive function $F$ on $\Gamma$ satisfying $F(\gamma)=0$ for $\gamma\in X$ and $F(\gamma)=1$ for $\gamma\notin\bigcup_{k=1}^K(U_2-\omega_{\nu_k})$.
Then the continuous function $G=F+\sum_{k=1}^K g_k$ on $\Gamma$ satisfies $G(\gamma)\geq 1$ for any $\gamma\in\Gamma$ since $F,g_1,g_2,\ldots,g_K$ are positive functions, and $F(\gamma)=1$ for $\gamma\notin\bigcup_{k=1}^K(U_2-\omega_{\nu_k})$, and $g_k(\gamma)=1$ for $\gamma\in U_2-\omega_{\nu_k}$.
Set $f_k=g_k/G$ for $k=1,2,\ldots,K$.
Then for $k=1,2,\ldots,K$, the positive function $f_k\in C_0(\Gamma)$ satisfies $f_k(\gamma)=0$ for any $\gamma\notin U-\omega_{\mu_k}$.
For $\gamma\in X$, we have 
\begin{align*}
\sum_{k=1}^K f_k(\gamma)&=\sum_{k=1}^K\frac{g_k(\gamma)}{G(\gamma)}\\
&=\frac{\sum_{k=1}^K g_k(\gamma)}{F(\gamma)+\sum_{k=1}^K g_k(\gamma)}\\
&=1.
\end{align*}

Since $X\subsetneqq\bigcup_{k=1}^K (U_2-\omega_{\nu_k})$, there exists $\gamma_0\notin X$ that is an element of $U_2-\omega_{\nu_{k_0}}$ for some $k_0\in\{1,2,\ldots,K\}$. 
Since $U_2-\omega_{\nu_{k_0}}$ is open and $X$ is closed,
we can choose an open set $O$ such that $\gamma_0\in O\subset U_2-\omega_{\nu_{k_0}}$ and $O\cap X=\emptyset$.
Let us take a positive function $f$ such that $f(\gamma)=0$ for any $\gamma\notin O$ and $f(\gamma_0)+\sum_{k=1}^K f_k(\gamma_0)$ is neither 0 nor 1.
Then $f_{k_0}'=f_{k_0}+f$ still satisfies that $f_{k_0}'(\gamma)=0$ for any $\gamma\notin U-\omega_{\mu_k}$.
We denote this new function $f_{k_0}'$ by $f_{k_0}$.
Then $K$ functions $f_1,f_2,\ldots,f_K$ satisfy $\sum_{k=1}^K f_k(\gamma)=1$ for $\gamma\in X$ and $\sum_{k=1}^K f_k(\gamma_0)\neq 0,1$. 
For $k=1,2,\ldots,K$, since $\sigma_{-\omega_{\mu_k}}f_k(\gamma)=0$ for any $\gamma\notin U\subset X$, the support of $\sigma_{-\omega_{\mu_k}}f_k$ is contained in $X$.
We get desired elements $f_1,f_2,\ldots,f_K\in C_0(\Gamma)$ and $\mu_1,\mu_2,\ldots,\mu_K\in\W_n$.
\eprf

\bpr\label{se}
If $\omega$ satisfies that $\Gamma=\overline{\{\omega_{\mu}\mid\mu\in\W_n\}}$,
then $\cp$ has a scaling element.
\epr

\bprf
Let $X$ be a compact neighborhood of $0\in\Gamma$ that differs from $\Gamma$. 
Let us take positive functions $f_1,f_2,\ldots,f_K\in C_0(\Gamma)$ and $\mu_1,\mu_2,\ldots,\mu_K\in\W_n$ that satisfy the four conditions in Lemma \ref{partition}.
Let us define $x=\sum_{k=1}^K S_{\mu_k}f_k^{1/2}\in\cp$.
Since $S_{\mu_k}^*S_{\mu_l}=\delta_{k,l}$,
$$x^*x=\sum_{k,l=1}^K f_k^{1/2}S_{\mu_k}^*S_{\mu_l}f_l^{1/2}=\sum_{k=1}^K f_k.$$
On the other hand, 
$$xx^*=\sum_{k,l=1}^K \big(S_{\mu_k}f_k^{1/2}f_l^{1/2}S_{\mu_l}^*\big)
=\sum_{k,l=1}^K\big((\sigma_{-\omega_{\mu_k}}f_k^{1/2})(\sigma_{-\omega_{\mu_k}}f_{l}^{1/2})S_{\mu_k}S_{\mu_l}^*\big).$$
Since the support of $\sigma_{-\omega_{\mu_k}}f_k^{1/2}$ is contained in $X$ for any $k=1,2,\ldots,K$ and $\sum_{k=1}^K f_k(\gamma)=1$ for $\gamma\in X$, we have $(x^*x)(xx^*)=xx^*$.

Finally we show $x^*x\neq xx^*$. If $x^*x=xx^*$, then $x^*x$ would become a projection. However, $x^*x=\sum_{k=1}^K f_k$ is not a projection, since there exists $\gamma_0\in\Gamma$ with $\sum_{k=1}^K f_k(\gamma_0)\neq 0,1$.
Thus $x$ is a scaling element.
\eprf

Since $\cp$ is simple, it has an infinite projection by Proposition \ref{bc1} 
and Proposition \ref{se}.
To prove that every non-zero hereditary subalgebra of $\cp$ has an infinite 
projection, we need the following lemma.
In the proof of it, we use some computations done in \cite{Ka} 
which is not difficult to see.
Let $\beta:\T\curvearrowright\cp$ be the gauge action defined by 
$\beta_t(S_\mu fS_\nu^*)=t^{|\mu|-|\nu|}S_\mu fS_\nu^*$,
and $E$ be the faithful conditional expectation of $\cp$ defined by 
$E(x)=\int_{\T}\beta_t(x)dt$ where $dt$ is the normalized Haar measure of $\T$.

\blem\label{pilem}
Let $y$ be a non-zero positive element of $\cp$, 
given as $y=\sum_{l=1}^LS_{\mu_l} f_lS_{\nu_l}^*$.
Let $C$ be a positive number with $1/\|E(y)\|< C^2$.
Then, there exist $a\in\cp$ with $\|a\|\leq C$ and an open set $O$ of $\Gamma$ 
such that $a^*ya$ becomes an element of $C_0(\Gamma)$ which is $1$ on $O$.
\elem

\bprf
Set $k=\max\{|\mu_l|,|\nu_l|\mid l=1,2,\ldots,L\}$ and 
$${\cal F}_k=\spa\{S_\mu fS_\nu^*\mid \mu,\nu\in\W_n^{(k)},\ f\in C_0(\Gamma)\}.$$
The $C^*$-algebra ${\cal F}_k$ is isomorphic to $C_0(\Gamma,\M_{n^k})$
and we will identify them.
We can see that $E(y)=\sum_{|\mu_l|=|\nu_l|}S_{\mu_l} f_lS_{\nu_l}^*$ and 
$E(y)\in {\cal F}_k$.
Set $u=\sum_{\mu\in\W_n^{(k)}}S_\mu S_1^kS_2S_\mu^*\in\On\subset M(\cp)$.
Routine computation shows that $u$ is an isometry and 
$u^*yu=\sigma_\gamma(E(y))$ where $\gamma=k\omega_1+\omega_2$.
Hence $u^*yu$ is a positive element of ${\cal F}_k$ whose norm is equal to 
$\|E(y)\|$.
One can find $\gamma_0\in\Gamma$ such that the norm of 
$(u^*yu)(\gamma_0)\in\M_{n^k}$ is $\|E(y)\|$.
The $C^*$-subalgebra $\spa\{S_\mu S_\nu^*\mid \mu,\nu\in\W_n^{(k)}\}$
of $\On\in M(\cp)$ is isomorphic to $\M_{n^k}$ and can be considered as 
the set of constant functions of $C_b(\Gamma,\M_{n^k})\cong M({\cal F}_k)$.
Take an element $\mu$ in $\W_n^{(k)}$ arbitrarily.
Then $S_{\mu}S_{\mu}^*\in M(\cp)$ is a minimal projection of $\M_{n^k}$.
Since $u^*yu$ is positive, 
$(u^*yu)(\gamma_0)$ is a positive element of $\M_{n^k}$.
Hence, there exists a partial isometry 
$v\in\spa\{S_\mu S_\nu^*\mid \mu,\nu\in\W_n^{(k)}\}$ such that 
$v^*v=S_{\mu}S_{\mu}^*$ and 
$$(v^*u^*yuv)(\gamma_0)=\|E(y)\|S_{\mu}S_{\mu}^*.$$
There exists a function $f\in C_0(\Gamma)$ with $v^*u^*yuv=S_{\mu}fS_{\mu}^*$,
because the projection $S_{\mu}S_{\mu}^*$ is minimal.
Since $f(\gamma_0)=\|E(y)\|$, there exists a positive function 
$g\in C_0(\Gamma)$ with $\|g\|\leq C$ such that 
$fg^2\in C_0(\Gamma)$ is 1 on some open neighborhood $O$ of $\gamma_0$.
If we set $a=u v S_{\mu}g\in\cp$, then, we get $\|a\|\leq C$ 
and $a^*ya=gfg$ becomes an element of $C_0(\Gamma)$ which is $1$ on $O$.
\eprf

\bthm\label{pi}
If $\omega$ satisfies that $\Gamma=\overline{\{\omega_{\mu}\mid\mu\in\W_n\}}$, then $\cp$ is simple and purely infinite.
\ethm

\bprf
To prove that $\cp$ is purely infinite, it suffices to show that 
there exists an infinite projection in the hereditary subalgebra 
$\overline{x(\cp)x}$ generated by $x$ for any non-zero positive element 
$x\in\cp$.
Take a non-zero positive element $x\in\cp$ 
and a sufficiently small positive number $\e>0$.
There exists a positive element $y$ with $\|x-y\|<\e$ 
that is a linear combination of elements of the form $S_\mu fS_\nu^*$.
Since $\|E(x)-E(y)\|\leq\|x-y\|<\e$, there exists a real number $C$ with 
$1/\|E(y)\|< C^2$ which depends only on $x$.
By Lemma \ref{pilem}, there exist $a\in\cp$ with $\|a\|\leq C$ 
and an open set $O$ of $\Gamma$ such that $a^*ya$ becomes an element of 
$C_0(\Gamma)$ which is $1$ on $O$.
Take an open subset $O_1$ of $O$ and a neighborhood $O_2$ of $0\in\Gamma$ 
with $O_1+O_2\subset O$.
Let $h$ be a non-zero positive function of $C_0(\Gamma)$ whose support is 
contained in $O_1$.
The crossed product $\cp$ has an infinite projection $p$
by Proposition \ref{bc1} and Proposition \ref{se}.
Since $\cp$ is simple and $p$ is a projection, there exist 
$x_1,x_2,\ldots,x_K\in\cp$ satisfying $\sum_{k=1}^K x_k^*h x_k=p$ 
by Lemma \ref{d}.
By Lemma \ref{type-}, we can choose $\mu_1,\mu_2,\ldots,\mu_K\in\W_n$ 
such that $S_{\mu_k}^*S_{\mu_l}=\delta_{k,l}$ 
and $\omega_{\mu_k}\in O_2$ for $k=1,2,\ldots,K$.
Set $b=\sum_{k=1}^K S_{\mu_k}h^{1/2} x_k$. We have
$$b^*b=\sum_{k,l=1}^K x_k^* h^{\frac12} S_{\mu_k}^*S_{\mu_l}h^{\frac12} x_l
=\sum_{k=1}^K x_k^* h x_k=p.$$
Since the support of $\sigma_{-\omega_{\mu_k}}(h^{1/2})$ is contained in $O$ 
for $k=1,2,\ldots,K$, and the function $a^*ya\in C_0(\Gamma)$ is 1 on $O$, 
we have $(a^*ya)b=b$.
Therefore, we get $b^*a^*yab=p$.
Thus $q=(y^{1/2}ab)(b^*a^*y^{1/2})$ is an infinite projection 
because it is equivalent to the infinite projection $p$.
The hereditary subalgebra $\overline{x(\cp)x}$ has a positive element 
$c=x^{1/2}abb^*a^*x^{1/2}$ which is close to an infinite projection $q$.
If we choose $\e>0$ so small that $\|q-c\|<1/2$,
then we get a projection $q_0=\chi (c)$ in $\overline{x(\cp)x}$ 
by the functional calculus where $\chi$ is a characteristic function 
of a certain neighborhood of 1.
The projection $q_0$ of $\overline{x(\cp)x}$ is infinite 
since it is close to an infinite projection $q$.
Therefore, $\cp$ is purely infinite.
\eprf

Once noting that $\cp$ is simple if and only if 
the closed semigroup generated by $\omega_1,\omega_2,\ldots,\omega_n$ 
and $-\omega_i$ is equal to $\Gamma$ for any $i=1,2,\ldots,n$
(see \cite[Theorem 4.4]{Ki} or \cite[Theorem 4.8]{Ka}),
we have the following corollaries.

\bco\label{dichotomy}
The crossed product $\cp$ is either purely infinite or AF-embeddable 
when it is simple.
\eco

\bco\label{pureinf} 
The crossed product $\cp$ is simple and purely infinite if and only if 
$\Gamma=\overline{\{\omega_{\mu}\mid\mu\in\W_n\}}$.
\eco

\bre
When the group $G$ is compact, 
crossed products $\cp$ are graph algebras \cite{KP}.
From this fact, one can easily prove Proposition \ref{AF!} and 
two corollaries above when the group $G$ is compact 
(see \cite{BPRS}, for example).
\ere

\bre
When the group $G$ is discrete, crossed products $\cp$ are never AF-embeddable
and Corollary \ref{dichotomy} implies that crossed products $\cp$ is 
purely infinite if it is simple.
This fact was already proved in \cite[Lemma 10]{KK2}.
\ere

\section{AF-embeddability of $\cpi$}

In this section, we deal with crossed products of the Cuntz algebra $\Oi$ 
which is the universal $C^*$-algebra generated by infinitely many isometries 
$S_1,S_2,\ldots$ satisfying $S_i^*S_j=\delta_{i,j}$.
Let us denote by $\W_\infty$ the set of words whose letters are 
$\{1,2,\ldots\}$, which is naturally identified with 
$\bigcup_{n=2}^\infty \W_n$.
We can define an isometry $S_\mu\in\Oi$ for $\mu\in\W_\infty$.
As in the case of $\On$, we define the action $\alpha^\omega$
of abelian group $G$ on $\Oi$ by 
$$\alpha^\omega_t(S_i)=\ip{t}{\omega_i}S_i\quad (i=1,2,\ldots,\ t\in G)$$
for $\omega=(\omega_1,\omega_2,\ldots)\in\Gamma^\infty$.
The crossed product $\cpi$ has the $C^*$-algebra 
$\C 1{\rtimes_{\alpha^\omega}}G$ which is isomorphic to $C_0(\Gamma)$.
One can easily see that $fS_\mu =S_\mu\sigma_{\omega_\mu}f$ 
for any $f\in C_0(\Gamma)\subset\cpi$ and any $\mu\in\W_\infty$, and
$$\cpi=\cspa\{ S_\mu fS_\nu^*\mid \mu,\nu\in\W_\infty,\ f\in C_0(\Gamma)\}.$$

\bpr\label{embedinfty}
If $\omega\in\Gamma^\infty$ satisfies 
$-\omega_i\notin\overline{\{\omega_{\mu}\mid\mu\in\W_n\subset\W_\infty\}}$ 
for any $i$ and any $n\in\N$, 
then the crossed product $\cpi$ is AF-embeddable.
\epr

\bprf
Fix an open base $\{U_i\}_{i\in\I}$ such that for any $i\in\I$, 
$\overline{U_i}$ is compact and for any $i\in\I$ and $\mu\in\W_\infty$, 
there exists $j\in\I$ with $U_j=U_i-\omega_\mu$.
Let $D_0(\Gamma)$ be the $C^*$-algebra generated $\{\chi_{U_i}\}_{i\in\I}$
in $L^\infty(\Gamma)$ and define the $C^*$-subalgebra $A$ of 
$B(H\otimes L^2(\Gamma))$ by
$$A=\cspa\{S_\mu fS_\nu^*\mid \mu,\nu\in\W_\infty,\ f\in D_0(\Gamma)\}.$$
The crossed product $\cpi$ can be embedded into $A$.
For a positive integer $n$ and a finite set $\lambda\subset\I$, 
we denote by $A_{\lambda,n}$ the $C^*$-subalgebra of $A$ generated by 
$$\{S_\mu \chi_{U_i}S_\nu^*\mid \mu,\nu\in\W_n\subset\W_\infty,\ i\in\lambda\}.$$
One can easily see that $A=\varinjlim A_{\lambda,n}$.
Take a positive integer $n$ and a finite set $\lambda\subset\I$ and fix them.
Since 
$-\omega_i\notin\overline{\{\omega_{\mu}\mid\mu\in\W_n\subset\W_\infty\}}$ 
for any $i$, there exists $K\in\N$ such that $pS_\mu p=0$ 
for any $\mu\in\W_n\subset\W_\infty$ with $|\mu|>K$ by Lemma \ref{shift}, 
where $p$ is the characteristic function of $\bigcup_{i\in\lambda}U_i$.
Once fixing such an integer $K$, we can define the projection 
$q\in A_{\lambda,n}$ in the same manner as in Section \ref{AFE} and
prove the same statement as in Lemma \ref{divide1} and Lemma \ref{divide2}.
Hence as in a similar way to Proposition \ref{K}, we can prove that
$A_{\lambda,n}$ is isomorphic to a direct product of finitely many $\K$.
Hence $A$ is an AF-algebra.
Since the crossed product $\cpi$ can be embedded into $A$, 
it is AF-embeddable.
\eprf

In the case of $\On$, we have the dichotomy (Corollary \ref{dichotomy}).
However in the case of $\Oi$, instead of dichotomy we have the following.

\bpr
For $\omega\in\Gamma^\infty$, the following are equivalent:
\benu
\item $\Gamma=\overline{\{\omega_{\mu}\mid\mu\in\W_\infty\}}$.
\item $\cpi$ is simple.
\item $\cpi$ is simple and purely infinite.
\eenu
\epr

\bprf
The equivalence between (i) and (ii) was proved in \cite{Ki}.
Obviously (iii) implies (ii).
One can prove the implication (i) $\Rightarrow$ (iii) in a similar way 
to arguments in Section \ref{secpi}, 
though we need more complicated computations to prove the proposition 
corresponding to Lemma \ref{pilem}.
\eprf


\begin{thebibliography}{BPRS}

\bibitem[BPRS]{BPRS}
Bates, T.; Pask, D.; Raeburn, I.; Szyma\'nski, W. {\it The $C\sp *$-algebras of row-finite graphs.} New York J. Math. {\bf 6} (2000), 307--324.

\bibitem[BC]{BC}
Blackadar, B. E.; Cuntz, J. {\it The structure of stable algebraically simple $C\sp{*} $-algebras.} Amer. J. Math. {\bf 104} (1982), no. 4, 813--822.

\bibitem[B1]{Br1}
Brown, N. P. {\it AF embeddability of crossed products of AF algebras by the integers.} J. Funct. Anal. {\bf 160} (1998), no. 1, 150--175.

\bibitem[B2]{Br2}
Brown, N. P. {\it Crossed products of UHF algebras by some amenable groups.} Hokkaido Math. J. {\bf 29} (2000), no. 1, 201--211.

\bibitem[B3]{Br3}
Brown, N. P. {\it On quasidiagonal $C\sp *$-algebras.} Preprint.

\bibitem[D]{D}
Davidson, K. R. {\it $C\sp *$-algebras by example.}
Fields Institute Monographs, {\bf 6}. American Mathematical Society, Providence, RI, 1996.

\bibitem[E]{E}
Evans, D. E. {\it On $O\sb{n}$.} Publ. Res. Inst. Math. Sci. {\bf 16} (1980), no. 3, 915--927.

\bibitem[Ka]{Ka}
Katsura, T. {\it The ideal structures of crossed products of Cuntz algebras by quasi-free actions of abelian groups.} Preprint.

\bibitem[Ki]{Ki}
Kishimoto, A. {\it Simple crossed products of $C\sp{*} $-algebras by locally compact abelian groups.} Yokohama Math. J. {\bf 28} (1980), no. 1-2, 69--85.

\bibitem[KK1]{KK1}
Kishimoto, A.; Kumjian, A. {\it Simple stably projectionless $C\sp *$-algebras arising as crossed products.} Canad. J. Math. {\bf 48} (1996), no. 5, 980--996.

\bibitem[KK2]{KK2}
Kishimoto, A.; Kumjian, A. {\it Crossed products of Cuntz algebras by quasi-free automorphisms.} Operator algebras and their applications, 173--192, Fields Inst. Commun., {\bf 13}, Amer. Math. Soc., Providence, RI, 1997.

\bibitem[KP]{KP}
Kumjian, A.; Pask, D. {\it $C\sp *$-algebras of directed graphs and group actions.} Ergodic Theory Dynam. Systems {\bf 19} (1999), no. 6, 1503--1519.

\bibitem[G]{Gr}
Green, P. {\it The structure of imprimitivity algebras.} J. Funct. Anal. {\bf 36} (1980), no. 1, 88--104.

\bibitem[Pi1]{Pi1}
Pimsner, M. V. {\it Embedding some transformation group $C\sp{*} $-algebras into AF-algebras.} Ergodic Theory Dynam. Systems {\bf 3} (1983), no. 4, 613--626.

\bibitem[Pi2]{Pi2}
Pimsner, M. V. {\it Embedding covariance algebras of flows into AF-algebras.} Ergodic Theory Dynam. Systems {\bf 19} (1999), no. 3, 723--740. 

\bibitem[PV]{PV}
Pimsner, M.; Voiculescu, D. {\it Imbedding the irrational rotation $C\sp{*} $-algebra into an AF-algebra.} J. Operator Theory {\bf 4} (1980), no. 2, 201--210.

\bibitem[Pu]{Pu}
Putnam, I. F. {\it The $C\sp *$-algebras associated with minimal homeomorphisms of the Cantor set.} Pacific J. Math. {\bf 136} (1989), no. 2, 329--353.

\bibitem[R]{Ro}
R\o rdam, M. {\it A simple $C^*$-algebra with a finite and an infinite projection.} Preprint.

\bibitem[V]{V}
Voiculescu, D. {\it Almost inductive limit automorphisms and embeddings into AF-algebras.} Ergodic Theory Dynam. Systems {\bf 6} (1986), no. 3, 475--484.


\end{thebibliography}
\end{document}